\documentclass[12pt]{article}
\usepackage[margin=2cm]{geometry}
\usepackage{amsfonts}
\usepackage{amssymb}
\usepackage{amsthm}
\usepackage{amsmath}
\usepackage{latexsym}
\usepackage{color}
\usepackage{mathrsfs}

\begin{document}
\allowdisplaybreaks[4]
\newtheorem{thme}{Theorem}
\newtheorem{lemma}{Lemma}
\newtheorem{pron}{Proposition}
\newtheorem{re}{Remark}
\newtheorem{thm}{Theorem}
\newtheorem{Corol}{Corollary}
\newtheorem{exam}{Example}
\newtheorem{defin}{Definition}
\newtheorem{remark}{Remark}
\newtheorem{property}{Property}
\newcommand{\bco}{\color{blue}}
\newcommand{\rco}{\color{red}}
\newcommand{\la}{\frac{1}{\lambda}}
\newcommand{\sectemul}{\arabic{section}}
\renewcommand{\thethme}{\sectemul.\Alph{thme}}
\renewcommand{\theequation}{\sectemul.\arabic{equation}}
\renewcommand{\thepron}{\sectemul.\arabic{pron}}
\renewcommand{\thelemma}{\sectemul.\arabic{lemma}}
\renewcommand{\there}{\sectemul.\arabic{re}}
\renewcommand{\thethm}{\sectemul.\arabic{thm}}
\renewcommand{\theCorol}{\sectemul.\arabic{Corol}}
\renewcommand{\theexam}{\sectemul.\arabic{exam}}
\renewcommand{\thedefin}{\sectemul.\arabic{defin}}
\renewcommand{\theremark}{\sectemul.\arabic{remark}}
\def\REF#1{\par\hangindent\parindent\indent\llap{#1\enspace}\ignorespaces}
\def\lo{\left}
\def\ro{\right}
\def\be{\begin{equation}}
\def\ee{\end{equation}}
\def\beq{\begin{eqnarray*}}
\def\eeq{\end{eqnarray*}}
\def\bea{\begin{eqnarray}}
\def\eea{\end{eqnarray}}
\def\r{random walk}
\def\o{\overline}

\title{\large\bf On the long tail property of product convolution
\thanks{Research supported by the National Science Foundation of China
(No. 11071182).}}
\author {Zhaolei Cui$^{1)}$\qquad Guancheng Jiang$^{2)}$\qquad Yuebao Wang$^{3)}$ \thanks{Corresponding author. \qquad Telephone: +86 512 67422726. \qquad Fax: +86 512 65112637. \quad E-mail: ybwang@suda.edu.cn (Y. Wang)}
\qquad\\ {\small\it 1) School of mathematics and statistics, Changshu Institute of Technology, Suzhou, P. R. China, 215000}\\
{\small\it 2) School of Economics£¬Shanghai University of Finance and Economics, P. R. China, 200433}\\
{\small\it 3) School of Mathematical Sciences, Soochow University, Suzhou, P. R. China, 215006}\\
}
\date{}

\maketitle

\begin{center}
{\noindent\small {\bf Abstract }}
\end{center}

{\small
Let $X$ and $Y$ be two independent random variables with corresponding distributions $F$ and $G$ supported on $[0,\infty)$. The distribution of the product $XY$, which is called the product convolution of $F$ and $G$, is denoted by $H$.
In this paper, some suitable conditions about $F $ and $G $ are given, under which
the distribution $H$ belongs to the long-tailed distribution class.
Here, $F$ is a generalized long-tailed distribution and is not necessarily an exponential distribution. Finally, a series of examples are given to show that the above conditions are satisfied by many distributions and one of them is necessary in some sense.
\medskip


\medskip

}

\section{Introduction and main results}
\setcounter{thm}{0}\setcounter{Corol}{0}\setcounter{lemma}{0}\setcounter{pron}{0}\setcounter{equation}{0}
\setcounter{re}{0}\setcounter{exam}{0}\setcounter{property}{0}\setcounter{defin}{0}

In this paper, without a special statement, let $X$ and $Y$ be two independent random variables with corresponding distributions $F$ and $G$ supported on $[0,\infty)$.
Denote the distribution of the product $XY$ by $H$ and it is called the product convolution of $F$ and $G$,
as opposed to the usual sum convolution $F*G$ defined to be the distribution of the sum $X+Y$.

As everyone knows, the concept of product convolution plays an important role in applied probability. For example, the present value of a random future claim, which is one of most fundamental quantities in finance and insurance, is expressed as the product of the random claim amount and the corresponding stochastic present value factor. Thus, the study of the tail behavior of product convolution has immediate implications in finance and insurance. There, the product convolution is often required to belong to the subexponential distribution class, see Embrechts and Goldie \cite{EG1980}, Cline \cite{C1986}, Cline and Samorodnitsky \cite{CS1994}, Tang and Tsitsiashvili \cite{TT2003, TT2004}, Tang \cite{T2006}, Tang \cite{T2008}, Liu and Tang \cite{LT2010}, Samorodnitsky and Sun \cite{SS2016}, Xu et al. \cite{XCWC2018}, etc. For this, as Cline and Samorodnitsky \cite{CS1994} points out: ``The first one provides conditions for $H$ to be `long-tailed', including the closure results for $\mathcal{L}$ (see below for its definition).'' The results in this regard can be found in Theorem 2.2 of Cline and Samorodnitsky \cite{CS1994}, Theorem 1.1 of Tang \cite{T2008}, etc. This paper focuses on this issue. To this end, the concepts and notations of some related distribution classes are first introduced as follows.

For a distribution $V$,
denote its tail distribution by $\overline{V} = 1-V$.
Unless otherwise stated, all limiting relations are according to $x\to\infty$.
With $s(f,g):=\limsup f(x)/g(x)$ for two positive functions
$f(\cdot)$ and $g(\cdot)$, we write $f(x)=o\big(g(x)\big)$ if $s(f,g)=0$, $f(x)=O\big(g(x)\big)$ if $s(f,g)<\infty$, $f(x)\lesssim g(x)$ if $s(f,g)\leq 1$, $f(x)\asymp g(x)$ if $s(f,g)<\infty$ and $s(g,f)<\infty$, and
$f(x)\sim g(x)$ if $s(f,g)=s(g,f)=1$.

A distribution $V$ supported on $(-\infty,\infty)$ is said to belong to the exponential distribution class $\mathcal{L}(\gamma)$ for some $\gamma\ge0$ if
\begin{eqnarray*}
\overline V(x-t)\sim e^{\gamma t}\overline V(x)\ \text{for each}\ t\in(-\infty,\infty).
\end{eqnarray*}
When $\gamma>0$ and the distribution $V$ is lattice, the variables $x$ and $t$ above should be restricted to values of the lattice span. When $\gamma=0$, it reduces to the class of long-tailed distributions and write $\mathcal{L}:=\mathcal{L}(0)$. It is well known that every $V\in\mathcal{L}$ is heavy tailed in the sense of infinite exponential moments. 

A distribution $V$ supported on $[0,\infty)$ is said to belong to the convolution equivalent distribution class $\mathcal{S}(\gamma)$ for some $\gamma\ge0$ if, $V\in\mathcal{L}(\gamma),\ \int_{[0,\ \infty)} e^{\gamma y}V(dy)<\infty$ and
\begin{eqnarray*}
\overline{V^{*2}}(x):=\overline{V*V}(x)\sim2\int_{[0,\ \infty)} e^{\gamma y}V(dy)\overline V(x).
\end{eqnarray*}
Furthermore, a distribution $V$ supported on $(-\infty,\infty)$ is still said to be convolution equivalent if the distribution $V_+$ defined by
$$V_+(x)=V(x)\textbf{1}_{\{x\ge0\}}(x)\ \ \text{for all}\ x\in(-\infty,\infty)$$
is convolution equivalent, where $\textbf{1}(A)$ denotes the indicator  of an event $A$, which is equal to $1$ if $A$ occurs and  to $0$ otherwise. In particular, the class $\mathcal{S}(0)=:\mathcal{S}$ is called as subexponential distribution class. When $\gamma=0$ and $V$ is supported on $[0,\infty)$, there is no requirement for $V$ to belong to the long-tailed distribution class.

The class $\mathcal{S}$ was introduced by Chistyakov \cite{C1964}. On its systematic discussion and application, see, for example, Embrechts et al. \cite{EMK1997}, Asmussen \cite{A2000} and Foss et al. \cite{FKZ2013}. The classes $\mathcal{L}(\gamma)$ and $\mathcal{S}(\gamma)$ for some $\gamma>0$ was introduced by Chover et al. \cite{CNW1973a,CNW1973b}. For some of the recent work involving these classes, see, for example, Watanabe \cite{W2008}, Watanabe and Yamamuro \cite{WY2017} and Cui et al. \cite{COWW2018}.


Another well-known class of heavy-tailed distributions is the class $\mathcal{D}$ of dominated variation tail introduced by Feller \cite{F1971}. Recall that a distribution $V$ supported on $(-\infty,\infty)$ belongs to the class $\mathcal{D}$ if
$$\overline {V}(tx)=O\big(\overline V(x)\big)\ \text{for each}\ t>1.$$



The class $\bigcup_{\gamma\ge0}\mathcal{S}(\gamma)\bigcup\mathcal{D}$ are properly included in the following distribution class introduced by Kl\"{u}ppelberg \cite{K1990}. A distribution $V$ supported on $(-\infty,\infty)$ is said to belong to the generalized subexponential distribution class $\mathcal{OS}$, if \begin{eqnarray*}\limsup\overline{V^{*2}}(x)/\overline V(x)=:C^*(V)<\infty.\end{eqnarray*}

The class $\bigcup_{\gamma\ge0}\mathcal{L}(\gamma)\bigcup\mathcal{D}$ is properly included in the generalized long-tailed distribution class $\mathcal{OL}$ introduced by Simura and Watanabe \cite{SW2005} if
\begin{eqnarray*}
\limsup\overline{V}(x-t)/\overline V(x)=:\limsup C^*(V,t,x)=:C^*(V,t)<\infty\ \text{for each}\ t>0.
\end{eqnarray*}

Some new results involving the two classes can be found in Beke et al. \cite{BBS2015}, Xu et al. \cite{XSW2015} and Xu et al. \cite{XSWC2016}, Wang et al. \cite{WXCY2018}, etc. Accordingly, we give the following indicator:
\begin{eqnarray*}
C_*(V,t):=\liminf\overline{V}(x-t)/\overline V(x)\ \ \text{for each}\ t>0.
\end{eqnarray*}
Clearly, $C^*(V,\cdot)\ \big(\text{or}\ C_*(V,\cdot)\big):[0,\infty)\mapsto[1,\infty)$ is a non-decreasing function in $t>0$, and
$$1\le C_*(V,\cdot)\le C^*(V,\cdot)\le\infty.$$ Thus $\lim_{t\to0} C^*(V,t)\ \big(\text{or}\ C_*(V,t)\big)$ and $\lim_{t\to\infty} C^*(V,t)\ \big(\text{or}\ C_*(V,t)\big)$ exist, which may be infinity and are denoted by $C^*(V,0)\ \big(\text{or}\ C_*(V,0)\big)$ and $C^*(V,\infty)\ \big(\text{or}\ C_*(V,\infty)\big)$, respectively.

Using terminology in Bingham et al. \cite{BGT1987}, we know that $V\in\mathcal{OL}$ if and only if the function $v(\cdot)$ belongs to the $O$-regularly varying function class $\mathcal{OR}$, where
$$v(x):=\overline{V}\big(\ln(x)\big)\ \text{for all}\ x\ge0$$
such that
$$0<v_*(u):=\liminf v(ux)/v(x)\le\limsup v(ux)/v(x)=:v^*(u)<\infty\ \text{for all}\ u\ge1.$$

Now we return to the long tail property of product convolution. The pioneering work in this area is attributed to Theorem 2.2 of Cline and Samorodnitsky \cite{CS1994}, and later a mature result is obtained by Theorem 1.1 of Tang \cite{T2008}. Here, we introduce the latter in more detail. We denote by $D[F]$ the set of all positive discontinuities of the distribution $F$ and recall the following condition:\\
\\
\textbf{Condition A} $D[F]=\emptyset$, or $D[F]\not=\emptyset$ and
\begin{equation}\label{blxtj}
\overline G(x/d)-\overline G\big((x+1)/d\big)=o\big(\overline H(x)\big)\ \text{for all}\ d\in D[F].
\end{equation}
In general, it is not easy to verify the condition (\ref{blxtj}) since $H$ is usually unknown. For this reason, Corollary 1.1 of Tang \cite{T2008} gives the following three useful sufficient conditions in terms of  the known distributions $F$ and $G$ for (\ref{blxtj}):

$(A)$ there is some $h > 0$ such that $G(x, x + h]$ is eventually non-increasing in $x$;

$(B)$ $G\in\mathcal{L}$;

$(C)$ it holds for all $c > 0$ that $\overline{G}(cx)/\overline{H}(x)\to0$, which is further implied by either $(C_1)$ $\overline{G}(vx)=o\big(\overline{G}(x)\big)$ for some $v>1$ or $(C_2)$ $\overline{G}(vx)=o\big(\overline{F}(x)\big)$ for some $v>0$.\\
Among them, Condition (B) can be slightly weakened to $G^{*k}\in\mathcal{L}$ for some integer $k\ge2$ by Theorem 2.1 (1a) of Xu et al. \cite{XFW2015}, see Xu et al. \cite{XCWC2018} for details.

\begin{thme}\label{1b}
Assume that $F\in\mathcal{L}(\gamma)$ for some $\gamma\ge0$,
then $H\in\mathcal{L}$ if and only if the Condition A holds.
\end{thme}

Based on Theorem 2.1 of Cline and Samorodnitsky \cite{CS1994} and Theorem \ref{1b}, Xu et al. \cite{XCWC2018} obtains an equivalent condition of $H\in\mathcal{S}$ under the prerequisite that $F\in\mathcal{S}$. The prerequisite of Theorem \ref{1b} is that $F\in\mathcal{L}(\gamma)$ for some $\gamma\ge0$. In practice, there are many distributions, which are not exponential, while their product convolution is still long-tailed, see Example \ref{exam31} and Example \ref{exam32} below. In this way, the following interesting question arises naturally.\\
\\
\textbf{Problem 1.1} If the distribution of $F$ and $G$ is not necessarily exponential, under what conditions can the product convolution $H$ be long-tailed?\\

For this problem, there are some results for some special distributions $F$ and $G$, see Theorem 3.1 and Theorem 3.2 of Xu et al. \cite{XCWC2018}, among which some conditions still involve the unknown distribution $H$. This paper hopes to get a general conclusion for the generalized long-tailed distribution $F$, in which all conditions only involve the known distributions $F$ and $G$. Here we emphasize that the class of generalized long-tailed distributions has a very wide range.\\


\begin{thm}\label{thm4}
For some $\alpha>0$, assume that
\begin{eqnarray}\label{thm41}
\lim e^{\lambda_1 x^\alpha}\overline{F}(x)=\infty\ \ for\ some\ \lambda_1>0,\ \lim e^{\lambda_2 x^\alpha}\overline{F}(x)=0\ \ for\ some\ \lambda_2>0
\end{eqnarray}
and
\begin{eqnarray}\label{thm42}
\lim e^{\lambda_3x^\alpha}\overline{G}(x)=0\ \ for\ some\ \lambda_3>0.
\end{eqnarray}
If $F\in\mathcal{OL}$, then $H\in\mathcal{OL}$. Further, $H\in\mathcal{L}$ if
$C^*(F,0)=1$. Conversely, if $H\in\mathcal{L}$, then $C_*(F,0)=1$.
\end{thm}
\begin{re}\label{re101}
Here, if $0<\alpha<1$, then $F$ is heavy-tailed; otherwise, $F$ is light-tailed. $G$ can be light-tailed or heavy-tailed. The tail of $G$ can be lighter or heavier than the tail of $F$. And $F$ and $G$ do not have to be exponential distributions, see Remark \ref{re103} below for details.
\end{re}
\begin{re}\label{re102}
When $F\in\mathcal{OL}$, the condition (\ref{thm41}) for $\alpha=1$ is automatically satisfied. In fact, by Theorem 2.2.7 of Bingham et al. \cite{BGT1987}, $F\in\mathcal{OL}$ if and only if
$$f(x):=\overline{F}(\ln(x))=\exp\Big\{\eta(x)+\int_{(1,x]}\big(\xi(y)/y\big)dy\Big\}\ \ for\ all\ x>1,$$
where $\xi(\cdot)$ and $\eta(\cdot)$ are bounded on $[x_0,\infty)$ for some $x_0>1$, and measurable.
It is easy to find the condition
(\ref{thm41}) holds for each pair of $\lambda_1>1>\lambda_2>0$.
\end{re}
\begin{re}\label{re10}
The condition $C^*(F,0)=1$ is necessary in some sense, see Example \ref{exam33} below. In particular, when $F\in\mathcal{L}(\gamma)$ for some $\gamma\ge0$, the condition is also automatically satisfied.
\end{re}
\begin{re}\label{re103}
i) When $\alpha=1$, there are many light-tailed distributions belonging to $\mathcal{OS}\setminus\bigcup_{\beta\geq0}\mathcal{L}(\beta)$ or $\mathcal{OL}\setminus\big(\bigcup_{\beta\geq0}\mathcal{L}(\beta)\bigcup\mathcal{OS}\big)$ such that $C^*(F,0)=1$, and satisfying the condition (\ref{thm41}) by Remark \ref{re102}. See Example \ref{exam31} and Example \ref{exam32} $i)$ below.

ii) When $0<\alpha<1$, there are many heavy-tailed distributions belonging to $\mathcal{OL}\setminus\big(\mathcal{L}\bigcup\mathcal{OS}\big)$ satisfying the condition $C^*(F,0)=1$, see Example \ref{exam32} $ii)$ below. In addition, the condition (\ref{thm41}) holds for these distributions.

Subsequently, if the condition (\ref{thm42}) also holds for the distribution $G$, then by Theorem \ref{thm4}, we have $H\in\mathcal{L}$.
\end{re}
\begin{re}\label{re101} Although the condition (\ref{thm41}) holds, the distribution $F$ may not belong to the class $\mathcal{OL}$. For example, if $\overline{F}(x)\sim e^{-x^2}$ and $G=F$, then the conditions (\ref{thm41}) and (\ref{thm42}) with $\alpha=2$ and $\lambda_1>1>\lambda_2>0$ hold, while $\overline{H}(x)=o\big(\overline{H}(x-t)\big)$ for each $t>0$. Therefore $H\notin\mathcal{OL}$ and $C^*(F,0)=\infty$. In fact, it follows from Lemma \ref{lemma201} and Lemma \ref{lemma202} below that, for any $t>0$,
$$\overline{H}(x-t)\sim\int_{(b_1(x),\ x/b_2(x)]}\overline{F}\big((x-t)/y\big)F(dy)\sim \int_{(b_1(x),\ x/b_2(x)]}\overline{F}\big(x/y\big)e^{2xt/y}F(dy)\ge e^{2tb_2(x)}\overline{H}(x),$$
where $b_1(\cdot)$ and $b_2(\cdot)$ are two positive functions satisfying $b_i(x)\uparrow\infty$, $b_i(x)/x\downarrow0,i=1,2$ and $b_1(x)=o\big(x/b_2(x)\big)$.
\end{re}

Next, we characterize the properties of product convolution from different angles.
\begin{pron}\label{pron11}
i) If $\lim x^{\delta_1}\overline{F}(x)=\infty$ for some $\delta_1>0$, then
\begin{eqnarray}\label{lemma2021}
\lim x^{\delta_1}\overline{H}(x)=\infty,
\end{eqnarray}
which is equivalent to that $\overline{V}(x)=o\big(\overline{H}(x)\big)$ for each distribution $V$ satisfying $\int_0^\infty y^\delta V(dy)<\infty$ for all $\delta>0$. At this time, $\overline{V}(cx)=o\big(\overline{H}(x)\big)$ for each $c>0$.
Further, if $F\in\mathcal{L}$ and $\int_0^\infty y^\delta G(dy)<\infty$ for all $\delta>0$, then $H\in\mathcal{L}$.

ii) For each $\delta>0$, if $\lim x^{\delta}\max\big\{\overline{G}(x),\overline{F}(x)\big\}=0$, then
\begin{eqnarray}\label{lemma2022}
\lim x^{\delta}\overline{H}(x)=0
\end{eqnarray}
and $H\notin\mathcal{D}$.
\end{pron}

In Section 2, we prove Theorem \ref{thm4} and Proposition \ref{pron11}. In order to illustrate the significance and value of Theorem \ref{thm4}, we give some examples of distributions satisfying the two conditions (\ref{thm41}) and (\ref{thm42}) and show that the condition $C^*(F,0)=1$ is necessary in some sense for $H\in\mathcal{L}$ in Section 3.

\section{Proofs of Theorem \ref{thm4} and Proposition \ref{pron11}}
\setcounter{thm}{0}\setcounter{Corol}{0}\setcounter{lemma}{0}\setcounter{pron}{0}\setcounter{equation}{0}
\setcounter{re}{0}\setcounter{exam}{0}\setcounter{property}{0}\setcounter{defin}{0}

\textbf{Proof of Theorem \ref{thm4}} We prove the theorem through the following two lemmas.

\begin{lemma}\label{lemma201} For some $\alpha>0$, assume that
\begin{eqnarray}\label{lemma201}
\lim e^{\lambda_1 x^\alpha}\overline{F}(x)=\infty\ \text{for some}\ \lambda_1>0,
\end{eqnarray}
then
\begin{eqnarray}\label{lemma202}
\lim e^{\lambda x^\alpha}\overline{H}(x)=\infty\ \text{for each}\ \lambda>0.
\end{eqnarray}
Thus, for each distribution $V$ satisfying
$\lim e^{\lambda_0x^\alpha}\overline{V}(x)=0$ for some $\lambda_0>0$,
we have
\begin{eqnarray}\label{lemma200}
\overline{V}(cx)=o\big(\overline{H}(x)\big)\ \text{for each}\ c>0.
\end{eqnarray}
\end{lemma}
\proof For any two constants $M>1$ and $0<\lambda<\lambda_1$, by (\ref{lemma201}), there is a positive number $x_1:=x_1(F,G,M,\lambda,\lambda_1,\alpha)$ large enough such that, for all $x\ge x_1$, it holds uniformly for all $a<y\le b$ that
$$e^{\lambda x^\alpha}\overline{F}(x/y)=e^{\lambda_1 x^\alpha/a^\alpha}\overline{F}(x/y) \ge e^{\lambda_1 (x/y)^\alpha}\overline{F}(x/y)\ge M/G[a,b),$$
where $a=(\lambda_1/\lambda)^{1/\alpha}$ and $b=2(\lambda_1/\lambda)^{1/\alpha}$. Thus
$$e^{\lambda x^\alpha}\overline{H}(x)\ge\int_{(a,b]}e^{\lambda_1(x/y)^\alpha}\overline{F}(x/y)G(dy)\ge M,$$
that is (\ref{lemma202}) holds for the constant $\lambda$.

For each $c>0$, we define a distribution $V_c$ such that $\overline{V_c}(x)=\overline{V}(cx)$ for all $x$. We take $\lambda=\lambda_0c^\alpha$, then
$$e^{\lambda x^\alpha}\overline{V_c}(x)=e^{\lambda_0(cx)^\alpha}\overline{V}(cx)\to0.$$
Therefore, by (\ref{lemma202}), we have
$$\overline{V}(cx)/\overline{H}(x)=e^{\lambda x^\alpha}\overline{V_c}(x)/e^{\lambda x^\alpha}\overline{H}(x)\to0,$$
that is (\ref{lemma200}) holds for the constant $c$.\hfill$\Box$


\begin{lemma}\label{lemma202} Assume that
\begin{eqnarray}\label{lemma203}
\max\big\{\overline{F}(cx),\ \overline{G}(cx)\big\}=o\big(\overline{H}(x)\big)\ \text{for each}\ c>0.
\end{eqnarray}
Then there are two functions $b_i(\cdot):[0,\infty)\mapsto(0,\infty)$ such that $b_i(x)\uparrow\infty$, $b_i(x)/x\downarrow0,i=1,2$, $b_1(x)=o\big(x/b_2(x)\big)$ and
\begin{eqnarray}\label{lemma204}
\overline{H}(x)\sim\int_{(b_1(x),\ x/b_2(x)]}\overline{F}(x/y)G(dy).
\end{eqnarray}
In addition, assume that $F\in\mathcal{OL}$, then it holds uniformly for all $0<t\le b_1(x)$ that
\begin{eqnarray}\label{lemma205}
\overline{H}\big(x-t\big)\sim\int_{(b_1(x),\ x/b_2(x)]}\overline{F}\big((x-t)/y\big)G(dy).
\end{eqnarray}
\end{lemma}
\proof Since $\overline{F}(cx)=o\big(\overline{H}(x)\big)$ for each $c>0$, for each integer $n\ge1$, there exists a sequence of positive numbers $\{x_0=1,x_n:=x_n(F):n\ge1\}$ such that $x_n\ge\max\{2x_{n-1},n^2\}$ and when $x\ge x_n$,
$$\overline{F}(x/n)<\overline{H}(x)/n.$$
Let $b_0(\cdot): [0,\infty)\mapsto(0,\infty)$ be a function such that
$$b_0(x)=\textbf{1}(0\le x <x_0=1)+\sum_{n=1}^\infty n\textbf{1}(x_n\le x < x_{n+1}).$$
Clearly, $b_0(x)\uparrow\infty,\ b_0(x)/x\rightarrow 0$ and $\overline{F}(x/b_0(x))=o\big(\overline{H}(x)\big)$. Define a continuous linear function $b_1(\cdot): [0,\infty)\mapsto(0,\infty)$ by \begin{eqnarray*}
b_1(x)&=&\textbf{1}(0\le x <x_0)+\sum_{n=1}^\infty \Big(x/(x_{n+1}-x_n)+\big(n-x_{n+1}/(x_{n+1}-x_n)\big)\Big)\textbf{1}(x_n\le x < x_{n+1})\nonumber\\ &=&\textbf{1}(0\le x <x_0)+\sum_{n=1}^\infty \big(n-(x_{n+1}-x)/(x_{n+1}-x_n)\big)\textbf{1}(x_n\le x < x_{n+1}). \end{eqnarray*}
Then for all $x\in[0,\infty)$, $b_0(x)\ge b_1(x)\uparrow\infty$. When $x_n\le x\le x_{n+1}$, $$b_1(x)/x=1/(x_{n+1}-x_n)+\big(n-x_{n+1}/(x_{n+1}-x_n)\big)\big/x\downarrow0.$$
And
$$\overline{F}\big(x/b_1(x)\big)\le\overline{F}\big(x/b_0(x)\big)=o\big(\overline{H}(x)\big).$$
Similarly, there exists a function $b_2(\cdot):[0,\infty)\mapsto(0,\infty)$ such that $b_2(x)\uparrow\infty$, $b_2(x)/x\downarrow0$, $x/b_2(x)>b_1(x)$ and
\begin{eqnarray}\label{lemma206}
\max\big\{\overline{F}(x/b_1(x)),\ \overline{G}(x/b_2(x))\big\}=o\big(\overline{H}(x)\big).
\end{eqnarray}
By (\ref{lemma206}), we have
\begin{eqnarray}\label{lemma207}
&&\max\Big\{\int_{[0,\ b_1(x)]}\overline{F}\big(x/y)G(dy),\int_{(x/b_2(x),\ \infty)}\overline{F}(x/y)G(dy)\Big\}\nonumber\\
&\le&\max\Big\{\overline{F}\big(x/b_1(x)\big),\overline{G}\big(x/b_2(x)\big)\Big\}=o\big(\overline{H}(x)\big).
\end{eqnarray}
Therefore, (\ref{lemma204}) holds.

For the above function $b_1(\cdot)$, by $F\in\mathcal{OL}$ and (\ref{lemma207}), it holds uniformly for all $0<t\le b_1(x)$ that
\begin{eqnarray*}
\int_{(b_1(x-t),\ b_1(x)]}\overline{F}\big((x-t)/y\big)G(dy)
&\le&\overline{F}\big(x/b_1(x)-1\big)\overline{G}\big(b_1(x-b_1(x))\big)\nonumber\\
&=&o\big(\overline{F}(x/b_1(x))\big)
=o\big(\overline{H}(x)\big).
\end{eqnarray*}
Therefore, combined with (\ref{lemma204}), (\ref{lemma205}) holds uniformly for all $0<t\le b_1(x)$.\hfill$\Box$

Now, we prove the theorem. Firstly, according to Lemma \ref{lemma201}, by (\ref{thm41}) and (\ref{thm42}), we have (\ref{lemma203}) holds. Then according to Lemma \ref{lemma202}, by $F\in\mathcal{OL}$, for each $t>0$ and $i=1,2$, there is a function $b_i(\cdot):[0,\infty)\mapsto(0,\infty)$ such that $b_i(x)\uparrow\infty$, $b_i(x)/x\downarrow0$, $b_1(x)=o\big(x/b_2(x)\big)$ and for any $\varepsilon>0$ and large enough $x$,
\begin{eqnarray*}
\overline{H}(x)&\le&\overline{H}(x-t)\sim\int_{(b_1(x),\ x/b_2(x)]}\overline{F}(x/y-t/y)G(dy)\nonumber\\
&\lesssim&\int_{(b_1(x),\ x/b_2(x)]}\overline{F}(x/y)C^*(F,t/y)G(dy)\nonumber\\
&\lesssim&C^*(F,\varepsilon)\overline{H}(x)\downarrow C^*(F,0)\overline{H}(x)\ \ \text{as}\ \varepsilon\downarrow0.
\end{eqnarray*}
Therefore, it is obvious that $H\in\mathcal{OL}$, $1\le C^*(H,t)\le C^*(F,0)$ for the constant $t>0$. Further, if $C^*(F,0)=1$, we have $H\in\mathcal{L}$.

Conversely, according to the uniform convergence theorem for the function class $\mathcal{OR}$, see, for example, Theorem 2.0.8 of Bingham et al. \cite{BGT1987}, by $F\in\mathcal{OL}$, for any $1>\varepsilon>0$ and large enough $x$, we have
\begin{eqnarray*}
\overline{H}(x-t)&\sim&\int_{(b_1(x),\ x/b_2(x)]}\overline{F}(x/y)
\Big(\overline{F}(x/y-t/y)\Big/\overline{F}(x/y)\Big)G(dy)\nonumber\\
&\ge&(1-\varepsilon)\int_{(b_1(x),\ x/b_2(x)]}\overline{F}(x/y)C_*(F,t/y)G(dy)\nonumber\\
&\gtrsim&(1-\varepsilon)C_*(F,0)\overline{H}(x)\uparrow C_*(F,0)\overline{H}(x)\ \ \text{as}\ \varepsilon\downarrow0.
\end{eqnarray*}
Combining with $H\in\mathcal{L}$ and $C_*(F,0)\ge1$, we know that $C_*(F,0)=1$.
\hfill$\Box$
\\
\\
\textbf{Proof of proposition \ref{pron11}}  i) If $x^{\delta_1}\overline{F}(x)\to\infty$ for some $\delta_1>0$, then
$$x^{\delta_1}\overline{H}(x)\ge x^{\delta_1}\int_{(x,\ \infty)}\overline{G}(x/y)F(dy)\ge \overline{G}(1)x^{\delta_1}\overline{F}(x)\to\infty.$$
Therefore, by Theorem (2) of Wang et al. \cite{WCY2005}, we know that (\ref{lemma2021}) is equivalent to the next conclusion.

For any $c>0$, We define a distribution $G_c$ as in the proof of Lemma \ref{lemma201}, then for any $\delta>0$,
$$\int_{[0,\ \infty)} y^\delta G_c(dy)=\int_{[0,\ \infty)} y^\delta dG(cy)=\int_{[0,\ \infty)}z^\delta G(dz)\Big/c^{\delta}<\infty.$$
Thus, by above conclusion, $\overline{G}(cx)=\overline{G_c}(x)=o\big(\overline{H}(x)\big)$. Further, if $F\in\mathcal{L}$, then $H\in\mathcal{L}$.


ii) For each $\delta>0$, since $x^\delta\max\{\overline{G}(x),\overline{F}(x)\}\to0$, $\sup_{x\ge0}x^\delta\overline{F}(x)<\infty$ and $EY^\delta<\infty$. Thus,
$$x^\delta\overline{H}(x)=\int_{[0,\ \infty)} y^\delta(x/y)^\delta\overline{F}(x/y)G(dy)<\infty$$
for the constant $\delta$. Therefore, (\ref{lemma2022}) holds by the arbitrariness of $\delta$, and $H\notin\mathcal{D}$.
\hfill$\Box$

\section{Some examples}

In this section, we first show that the necessity of the condition $C^*(F,0)=1$ in some sense. Next, we provide some examples of distribution $F$ in the class $\mathcal{OS}\setminus\bigcup_{\gamma\ge0}\mathcal{L}(\gamma)$ or $\mathcal{OL}\setminus\big(\bigcup_{\gamma\ge0}\mathcal{L}(\gamma)\bigcup\mathcal{OS}\big)$ with normal shapes satisfying the conditions (\ref{thm41}) and $C^*(F,0)=1$. It is easy to see that there are many such distributions. 
\begin{exam}\label{exam33}
Assume that $F_0\in\mathcal{L}(\gamma)$ for some $\gamma\ge0$ is a continuous distribution.
For example, we can take
\begin{eqnarray}\label{exam301}
\overline{F}_0(x)=\textbf{\emph{1}}(x<0)+e^{-\gamma x-x^{1/2}}\textbf{\emph{1}}(x\ge0)\ \ \ \text{for all}\  x\in(-\infty,\infty).
\end{eqnarray}
In fact, this distribution belongs to the class $\mathcal{S}(\gamma)$. We construct a distribution $F$ generated by $F_0$ as following:
\begin{eqnarray}\label{exam302}
\overline{F}(x)=\overline{F}_0(x)\textbf{\emph{1}}(x<x_1)+
\sum_{i=1}^{\infty}\Big(\overline{F}_0(x_i)\textbf{\emph{1}}(x_i\le x<y_i)+\overline{F}_0(x)\textbf{\emph{1}}(y_i\le x<x_{i+1})\Big)
\end{eqnarray}
for all $x\in(-\infty,\infty)$, where $\{y_i>x_i>0:i\ge1\}$ is a sequence of positive numbers
such that
$$\lim_{i\to\infty}\overline{F}_0(x_i)/\overline{F}_0(y_i)=a\ \ \text{and}\ \ x_{i+1}-y_{i}\ge x_{i}-y_{i-1}$$
for some $a>1$ and all $i\ge1$.

First, it is easy to find that, $\overline{F}(x)\asymp\overline{F}_0(x)$, $F\in\mathcal{OL}\setminus\bigcup_{\beta\ge0}\mathcal{L}(\beta)$, and the condition (\ref{thm41}) holds for any $\lambda_1>\gamma>\lambda_2>0$ in the case that $\gamma>0$. When $\gamma=0$, the condition (\ref{thm41}) also may be established. For example, the distribution $F$ generated by the above specific distribution $F_0$ in (\ref{exam301}), satisfies the condition (\ref{thm41}) for any $\lambda_1>1>\lambda_2>0$ with $\alpha=1/2$.

Next, we note that, when $\gamma>0$, there is a constant $b>0$ such that $y_i-x_i<b$ for all $i\ge1$. Otherwise, for $i$ large enough and any constant $(\ln6a)/\gamma<M\le y_i-x_i$,
$$1/(2a)<\overline{F}_0(x_i+y_i-x_i)/\overline{F}_0(x_i)\le\overline{F}_0(x_i+M)/\overline{F}_0(x_i)\le2e^{-\gamma M}<1/(3a).$$
Here there is a contradiction. However, when $\gamma=0$, if there is a constant $b>0$ such that $y_i-x_i<b$ for $i$ large enough, then
$$\liminf_{i\to\infty}\overline{F}_0(x_i+y_i-x_i)/\overline{F}_0(x_i)\ge\lim_{i\to\infty}\overline{F}_0(x_i+b)/\overline{F}_0(x_i)=1,$$
which is inconsistent with the fact that $\lim_{i\to\infty}\overline{F}_0(y_i)/\overline{F}_0(x_i)=1/a<1$, thus $y_i-x_i\to\infty$ as $i\to\infty$.

Third, we prove $C^*(F,0)=a$ and $C_*(F,0)=1$. For the former, we only need to calculate $C^*(F,t)$ for $t>0$ small enough such that, for all $i\ge1$, $\widehat{x_i}$ and $\widehat{x_i}-t$ are located in adjacent intervals for any $\widehat{x_i}>0$, and $y_i$ (or $x_i$) and $y_i+t$ (or $x_i-t$) are located in the same interval. When $\widehat{x_i}\in[y_i,x_{i+1})$ and $\widehat{x_i}-t\in[x_i,y_i)$,
\begin{eqnarray*}
\limsup_{i\to\infty}\overline{F}(\widehat{x_i}-t)/\overline{F}(\widehat{x_i})&=&\limsup_{i\to\infty}\big(\overline{F}_0(x_i)/\overline{F}_0(y_i)\big)
\big(\overline{F}_0(y_i)/\overline{F}_0(\widehat{x_i})\big)\nonumber\\
&\le&\lim_{i\to\infty} \big(\overline{F}_0(x_i)/\overline{F}_0(y_i)\big)
\big(\overline{F}_0(y_i)/\overline{F}_0(y_i+t)\big)\nonumber\\
&=&ae^{\gamma t}\nonumber\\
&=&a\lim_{i\to\infty} \overline{F}(y_i+t-t)/\overline{F}(y_i+t);
\end{eqnarray*}
when $\widehat{x_i}\in[x_i,y_i)$ and $\widehat{x_i}-t\in[y_{i-1},x_{i})$, for $i$ large enough,
$$\overline{F}(\widehat{x_i}-t)/\overline{F}(\widehat{x_i})=\overline{F}_0(x_i-t+\widehat{x_i}-x_i)
/\overline{F}_0(x_i)\le \overline{F}_0(x_i-t)/\overline{F}_0(x_i)\to e^{\gamma t}\ \text{as}\ i\to\infty.$$
By the above fact, we have
$$C^*(F,0)=\lim_{t\to0}C^*(F,t)=\lim_{t\to0}\lim_{i\to\infty}a\overline{F}(y_i)/\overline{F}(y_i+t)
=\lim_{t\to0}ae^{\gamma t}=a.$$

On the other hand, by $\lim_{i\to\infty}\overline{F}(x_i+t)/\overline{F}(x_i)=1$ for each $0<t<y_i-x_i$, we know that $C_*(F,0)=1$.
Further, we assume $G$ satisfies the condition (\ref{thm42}), then $H\in\mathcal{OL}$ by Theorem \ref{thm4}.

Finally, we show that the condition $C^*(F,0)=1$ of Theorem \ref{thm4} is necessary in some sense. In fact, there is a distribution $F$, which satisfies the above requirements, but its product convolution $H$ with $G=F$ does not belong to the class $\mathcal{L}$.

We take a specific distribution $F$ in (\ref{exam302}) with some $\gamma>0$ and $\{x_i=2i,y_i=2i+1:i\ge1\}$, which is generated by distribution $F_0$ in (\ref{exam301}), such that
$$\overline{F}(x)=\overline{F}_0(x)\textbf{\emph{1}}(x<2)+
\sum_{i=1}^{\infty}\Big(\overline{F}_0(2i)\textbf{\emph{1}}(2i\le x<2i+1)+\overline{F}_0(x)\textbf{\emph{1}}(2i+1\le x<2(i+1))\Big)$$
for all $x\in(-\infty,\infty)$. Clearly,
$$C^*(F,0)=a=\lim_{i\to\infty}\overline{F}_0(2i)/\overline{F}_0(2i+1)=e^{\gamma}>1$$
and $C_*(F,0)=1$. Thus $F\in\mathcal{OL}\setminus\bigcup_{\beta\ge0}\mathcal{L}(\beta)$ and the conditions (\ref{thm41}) and (\ref{thm42}) hold with $\alpha=1$ for each pair of $\lambda_1>\gamma>\lambda_2=\lambda_3>0$.

Take $G=F$, then by (\ref{lemma205}) of Lemma \ref{lemma202}, there exist two functions $b_i(\cdot):[0,\infty)\mapsto(0,\infty)$ satisfying $b_i(x)\uparrow\infty$, $b_i(x)/x\downarrow0$ , $i=1,2$ and $b_1(x)=o\big(x/b_2(x)\big)$, such that (\ref{lemma204}) holds.
Moreover, we denote $K_2(x)=\lfloor x/b_2(x)\rfloor$ and might as well assume that
it is odd, where the symbol $\lfloor x\rfloor$ stands for the largest integer not exceeding a real number $x$. Let $X$ be a r.v. with the distribution $F$. Then for each $t>0$, by (\ref{lemma204}), it still holds that
\begin{eqnarray}\label{exam303}
&&\overline{H}\big(x-t\big)\sim\sum_{b_1(x)<2i+1\le K_2(x)}\overline{F}\Big(\frac{x-t}{2i+1}\Big)P(X=2i+1)\nonumber\\
&&\ \ \ \ \ \ \ \ \ \ \ \ \ \ \ \ \ \ \ \ \ \ \ +\sum_{b_1(x)< 2i+1<2i+2<K_2(x)}\int_{(2i+1,2i+2)}\overline{F}\Big(\frac{x-t}{y}\Big)F_0(dy)\nonumber\\
&\ge&\sum_{b_1(x)< 2i+1\le K_2(x)}
\overline{F}\Big(\frac{x-t}{2i+1}\Big)P(X=2i+1)+\sum_{b_1(x)< 2i+1<2i+2<K_2(x)}\int_{(2i+1,2i+2)}\overline{F}(x/y)F_0(dy)\nonumber\\
&=:&H_1(x,t)+H_2(x,0).
\end{eqnarray}
Take a sequence$\{\widehat{x}_n=(2n+1)!!:n\ge1\}$, so we can know the number $\widehat{x}_n/(2i+1)$ also is odd, subsequently, it holds uniformly for all $b_1(x)<2i+1\le K_2(x)$ that
\begin{eqnarray}\label{exam304}
\overline{F}\Big(\frac{\widehat{x}_n-t}{2i+1}\Big)P(X=2i+1)\sim e^\gamma\overline{F}\Big(\frac{\widehat{x}_n}{2i+1}\Big)P(X=2i+1)\ \ \text{as}\ \ n\to\infty.
\end{eqnarray}
Similarly, it holds uniformly for all $b_1(x)<2i+1<2i+2\le K_2(x)$ that
\begin{eqnarray}\label{exam305}
\int_{(2i+1,2i+2)}\overline{F}(\widehat{x}_n/y)F_0(dy)&\leq& \overline{F}\big(\widehat{x}_n/(2i+3)\big)\big(\overline{F}_0(2i+1)-\overline{F}_0(2i+2)\big)\nonumber\\
&\sim&e^\gamma\overline{F}\big(\widehat{x}_n/(2i+3)\big)\big(\overline{F}_0(2i+2)-\overline{F}_0(2i+3)\big)\nonumber\\
&=& e^\gamma\overline{F}\big(\widehat{x}_n/(2i+3)\big)P(X=2i+3).
\end{eqnarray}
By (\ref{exam303}), (\ref{exam304}) and (\ref{exam305}), we have
\begin{eqnarray*}
\overline{H}(\widehat{x}_n-t)&\gtrsim&e^\gamma H_1(\widehat{x}_n,0)+H_2(\widehat{x}_n,0)\nonumber\\
&=&(e^\gamma-1)H_1(\widehat{x}_n,0)+H_1(\widehat{x}_n,0)+H_2(\widehat{x}_n,0)\nonumber\\
&\gtrsim&(e^\gamma-1)\big(H_1(\widehat{x}_n,0)+H_2(\widehat{x}_n,0)\big)/2+H_1(\widehat{x}_n,0)+H_2(\widehat{x}_n,0)\nonumber\\
&\sim&\big(1+(e^\gamma-1)/2\big)\overline{H}(\widehat{x}_n),
\end{eqnarray*}
which implies $H\notin\mathcal{L}$.
\end{exam}

\begin{exam}\label{exam31}
Here, we give a class of distributions, in which each distribution $F$ belongs to the class $\mathcal{OS}\setminus\bigcup_{\beta\ge0}\mathcal{L}(\beta)$, and satisfy the conditions (\ref{thm41}) and $C^*(F,0)=1$ with each pair of $(\lambda_1,\lambda_2)$ such that $\lambda_1>\gamma>\lambda_2>0$ for some $\gamma>0$.

Let $F_0$ be an absolutely continuous distribution on $[0,\infty)$ belonging to the class $\mathcal{L}(\gamma)\bigcap\mathcal{OS}$ for some $\gamma>0$ with a density $f_0$. For example, $$\overline{F_0}(x)=\textbf{\emph{1}}_{\{x<0\}}(x)+\textbf{\emph{1}}_{\{x\ge0\}}(x)/e^{\gamma x}(1+x)^3\ \  \text{for all}\ x.$$
By Karamata's theorem, we have $f_0(x)\sim \gamma\overline{F_0}(x)$. So we might as well assume that $$b:=\inf_{x\ge0}f_0(x)/\overline{F_0}(x)>0.$$
We take a constant $a>(1+b)/b$. Let $g(\cdot):=g(\cdot,F_0,a)$ be a function such that
$$g(x)=\overline{F_0}(x)(1+\sin x/a)\textbf{\emph{1}}_{\{x\ge0\}}(x)\ \ \text{for all}\ \ x\ge0.$$
Because $g(0)=1,\ g(\infty)=0$ and for $x\ge0$,
$$\frac{d}{dx}g(x)=-\overline{F_0}(x)\Big(\frac{f_0(x)}{\overline{F_0}(x)}\big(1+\frac{\sin x}{a}\big)+\frac{\cos x}{a}\Big)<-\overline{F_0}(x)\Big(b\big(1+\frac{\sin x}{a}\big)+\frac{\cos x}{a}\Big)<0,$$
the function $F(\cdot):=1-g(\cdot):=F(\cdot,F_0,a)$ is a distribution function supported on $[0,\infty)$.

Clearly, $F\in\mathcal{OS}\subset\mathcal{OL}$, while $F\notin\bigcup_{\beta\ge0}\mathcal{L}(\beta)$. Thus, for each $t>0$, there exists a sequence of positive numbers $\{x_n:n\ge1\}$ such that $x_n\uparrow\infty$ as $n\to\infty$ and
\begin{eqnarray*}
&&C^*(F,t)=\lim_{n\to\infty}\overline{F}(x_n-t)/\overline{F}(x_n)\nonumber\\
&=&\lim_{n\to\infty}e^{\gamma t}\big(a+\sin(x_n-\lfloor x_n/2\pi\rfloor2\pi-t)\big)\big/(a+\sin(x_n-\lfloor x_n/2\pi\rfloor2\pi))<\infty.
\end{eqnarray*}
Since $x_n-\lfloor x_n/2\pi\rfloor2\pi\in[0,2\pi)$ for $n\ge1$, there exists a subsequence $\{y_m:m\ge1\}$ of the sequence $\{x_n:n\ge1\}$ and a constant $c:=c(F,t)\in[0,2\pi]$ such that $y_m\uparrow\infty$ and $y_m-\lfloor y_m/2\pi\rfloor2\pi\to c$ as $m\to\infty$. Thus,
\begin{eqnarray*}
1<C^*(F,t)=e^{\gamma t}\big(a+\sin(c-t)\big)\big/(a+\sin c)\to1\ \ \text{as}\ t\to0.
\end{eqnarray*}

In addition, it is also obvious that, for each pair of $\lambda_1>\gamma>\lambda_2>0$, (\ref{thm41}) holds.\hfill$\Box$
\end{exam}

\begin{exam}\label{exam32}
There is a light-tailed distribution and a heavy-tailed distribution, both of which belong to the class $\mathcal{OL}\setminus\big(\bigcup_{\beta\ge0}\mathcal{L}(\beta)\bigcup\mathcal{OS}\big)$ and satisfy the condition (\ref{thm41}) for some $\lambda_1>\lambda_2>0$ and the condition $C^*(\cdot,0)=1$. 

Let $\theta\in
\big(3/2,(\sqrt{5}+1)/2\big)$ and $r=(\theta+1)/\theta$ be constants.
Assume $a>1$ is so large that $a^{r} > 8a$.
We define a distribution $F_0$ supported on $[0,\infty)$ such that
\begin{eqnarray}\label{exam01}
\overline{F_0}(x)&=&\textbf{\emph{1}}(x<a_0)+C\sum\limits_{n=0}^{\infty}\Big(\Big(\sum\limits_{i=n}^{\infty}
a_i^{-\theta}-a_n^{-\theta-1}(x-a_n)\Big)
\textbf{\emph{1}}_{\{x\in[a_n,2a_n)\}}(x)\nonumber\\
& &+\sum\limits_{i=n+1}^{\infty}a_i^{-\alpha}\textbf{\emph{1}}_{\{x\in[2a_n,a_{n+1})\}}(x)\Big),
\end{eqnarray}
where $C$ is a regularization constant and $a_n=a^{r^n}$ for all nonnegative integers. Here, $F_0$ is heavy-tailed and belongs to the class $\mathcal{OS}\setminus\mathcal{L}$, which comes from Proposition 5.1 of Xu et al. \cite{XWCY2018}.

In addition, we have
$C^*(F_0,t)=1+t$ for each $t>0$, thus $C^*(F_0,0)=1$. In fact, for any $\widehat{x}_n\in[a_n,2a_n)$, if $\widehat{x}_n-t$ is in the same interval, then by $a_{n+1}^{-\theta}=a_{n}^{-\theta-1}$, we have
\begin{eqnarray*}
\limsup\overline{F_0}(\widehat{x}_n-t)/\overline{F_0}(\widehat{x}_n)&=&\limsup_{n\to\infty}\lo(1+ta_n^{-\theta-1}
\Big/\Big(\sum\limits_{i=n}^{\infty}a_i^{-\theta}-a_n^{-\theta-1}(\widehat{x}_n-a_n)\Big)\ro)\nonumber\\
&\le&\limsup_{n\to\infty}\overline{F_0}(2a_n-t)/\overline{F_0}(2a_n)=1+t.
\end{eqnarray*}
For any $\widehat{x}_n\in[2a_n,a_{n+1})$, if $\widehat{x}_n-t\in[a_n,2a_n)$, then
\begin{eqnarray*}
\limsup\overline{F_0}(\widehat{x}_n-t)/\overline{F_0}(\widehat{x}_n)
&=&\limsup\overline{F_0}(\widehat{x}_n-t)/\overline{F_0}(2a_n)\nonumber\\
&=&\limsup_{n\to\infty}\Big(1+(t+a_n-\widehat{x}_n)a_n^{-\theta-1}\Big/\sum\limits_{i=n+1}^{\infty}a_i^{-\theta}\Big)\nonumber\\
&\le&\lim_{n\to\infty}\overline{F_0}(2a_n-t)/\overline{F_0}(2a_n)=1+t;
\end{eqnarray*}
if $\widehat{x}_n-t \in[2a_n,a_{n+1})$, then
$$\limsup\overline{F_0}(\widehat{x}_n-t)/\overline{F_0}(\widehat{x}_n)=1\le 1+t=\lim_{n\to\infty}\overline{F_0}(2a_n-t)/\overline{F_0}(2a_n).$$
By the above fact, we have
$$C^*(F_0,t)=\lim_{n\to\infty}\overline{F_0}(2a_n-t)/\overline{F_0}(2a_n)=1+t.$$

$i)$ For some $\gamma>0$, let $F$ be a distribution defined by
$$\overline{F}(x)=e^{-\gamma x}\overline{F_0}(x)\ \ \text{for all}\ \ x\in[0,\infty).$$
Then $F\in\mathcal{OL}\setminus\big(\bigcup_{\beta\ge0}\mathcal{L}(\beta)\bigcup\mathcal{OS})$ and
$$C^*(F,t)=(1+t)e^{\gamma t}\ \ \text{for each}\ \ t>0,$$
thus $C^*(F,0)=1$. Clearly, (\ref{thm41}) holds for each pair of $\lambda_1>\gamma>\lambda_2>0$ and $F$ is light-tailed.

$ii)$ For some $\gamma>0$, let $F$ be a distribution defined by
$$\overline{F}(x)=e^{-\gamma x^{1/2}}\overline{F_0}(x)\ \ \text{for all}\ \ x\in(-\infty,\infty).$$
Then $F$ is a heavy-tailed distribution belonging to the class $\mathcal{OL}\setminus\big(\mathcal{L}\bigcup\mathcal{OS})$ and
$$C^*(F,t)=1+t\ \ \text{for each}\ \ t>0,$$
thus $C^*(F,0)=1$. Clearly, (\ref{thm41}) holds for each pair of $\lambda_1>\gamma>\lambda_2>0$ with $\alpha=1/2$.
\hfill$\Box$
\end{exam}

\end{document}